\documentclass[reqno]{amsart}
\usepackage{amsmath, amsthm, amssymb, mathabx}
\begin{document}
\newtheorem{theo}{Theorem}[section]
\newtheorem{defin}[theo]{Definition}
\newtheorem{rem}[theo]{Remark}
\newtheorem{lem}[theo]{Lemma}
\newtheorem{cor}[theo]{Corollary}
\newtheorem{prop}[theo]{Proposition}
\newtheorem{exa}[theo]{Example}
\newtheorem{exas}[theo]{Examples}
%%%%%%%%%%%%%%%%%%%%%%%%%%%%%%%%%%%
%%%%%%%%%%%%%%%%%%%%%%%%%%%%%%%%%%%%
%
%
\thanks{}
\title[Radial and nonradial solutions of a strongly indefinite elliptic system on $\mathbb{R}^N$]{Radial and nonradial solutions of a strongly indefinite elliptic system on $\mathbb{R}^N$}

\author[C. J. Batkam]{Cyril Joel Batkam}
\address{Cyril Joel Batkam \newline
D\'epartement de math\'ematiques,
\newline
Universit\'e de Sherbrooke,
\newline
Sherbrooke, (Qu\'ebec), J1K 2R1, CANADA.}
\email{cyril.joel.batkam@usherbrooke.ca}

\maketitle
\begin{abstract}
This paper is concerned with the following system of elliptic equations
\begin{equation*}
    \left\{
                                                  \begin{array}{ll}
                                                    -\Delta u+u= F_u(|x|,u,v),  & \hbox{} \\ \\
                                                    -\Delta v+v=- F_v(|x|,u,v), & \hbox{} \\ \\
                                                    \,\,\,\,\,u,v\in H^1(\mathbb{R}^N). & \hbox{}
                                                  \end{array}
                                                \right.
\end{equation*}
It is shown that if $F$ is even in $(u,v)$ and satisfies some growth conditions, then the system has infinitely many both radial and nonradial solutions. The proof relies on the Principle of Symmetric Criticality and a generalized Fountain Theorem for strongly indefinite even functionals.
\end{abstract}

%
%%%%%%%%%%%%%%%%%%%%%%%%%%%%%%%%%%%%%%%%%%%%%%%%%%%%
%%%%%%%%%%%%%%%%%%%%%%%%%%%%%%%%%%%%%%%%%%%%%%%%%%%%
\section{Introduction}
%%%%%%%%%%%%%%%%%%%%%%%%%%%%%%%%%%%%%%%%%%%%%%%%%%%%
%%%%%%%%%%%%%%%%%%%%%%%%%%%%%%%%%%%%%%%%%%%%%%%%%%%%
 In this paper, we study the existence and multiplicity of solutions of the noncooperative elliptic system

\begin{equation*}
    (\mathcal{S})\,\,\,\,\,\,\,\,\,\,\left\{
                                                  \begin{array}{ll}
                                                    -\Delta u+u= F_u(x,u,v),  & \hbox{} \\ \\
                                                    -\Delta v+v=- F_v(x,u,v), & \hbox{} \\ \\
                                                    \,\,\,\,\,u,v\in H^1(\mathbb{R}^N), & \hbox{}
                                                  \end{array}
                                                \right.
\end{equation*}
where $F:\mathbb{R}^N\times \mathbb{R}^2\rightarrow \mathbb{R}$ is of class $\mathcal{C}^1$, and $F_w$ designates the partial derivative of $F$ with respect to $w$. The solutions of such a system are steady state of reaction-diffusion systems which modeled  many phenomena in Biology, in Chemical Reactions or in Physic.
\par It is well known that this problem has a variational structure, and its weak solutions are critical points of the functional $J$ defined on the Hilbert space $H^1(\mathbb{R}^N)\times H^1(\mathbb{R}^N)$ by
\begin{equation}\label{j2}
J(u,v):=\frac{1}{2}\int_{\mathbb{R}^N}\big(|\nabla u|^2+u^2\big)dx-\frac{1}{2}\int_{\mathbb{R}^N}\big(|\nabla v|^2+v^2\big)dx-\int_{\mathbb{R}^N} F(x,u,v)dx.
\end{equation}
We recall that the couple $(u,v)$ is a weak solution of $(\mathcal{S})$ if $(u,v)\in H^1(\mathbb{R}^N)\times H^1(\mathbb{R}^N)$, and satisfies for any $\phi,\varphi\in H^1(\mathbb{R}^N)$:
\begin{equation*}
    \int_{\mathbb{R}^N}\Big[\nabla u\nabla\phi+u\phi-\nabla v\nabla\varphi-u\varphi-\phi F_u(x,u,v)-\varphi F_v(x,u,v)\Big]dx=0.
\end{equation*}
The functional $J$ is strongly indefinite in the sense that it is neither bounded from below nor from above, even modulo subspaces of finite dimension or codimension. Therefore, the usual multiplicity critical point theorems such as the symmetric mountain pass theorem of Ambrosetti and Rabinowitz \cite{AR}, or the fountain  theorems of Bartsch and Willem \cite{W} are not applicable. This is the first difficulty to overcome when investigating the existence of solutions of $(\mathcal{S})$ by variational methods. The second difficulty is the lack of the compactness of the Sobolev embeddings, since we consider the whole space $\mathbb{R}^N$. Fortunately when the problem has some symmetry properties,
% for example when it is invariant by a group of orthogonal transformations,
  it suffices to consider invariant functions to recover compactness.
\par There are various methods in literature dealing with symmetric strongly indefinite functionals. Among others we can mention Rabinowitz \cite{Rabi}, Benci and Rabinowitz \cite{BenRabi}, Benci \cite{Ben}, Li \cite{Li}, Bartsch and Clapp \cite{BaCla}, Costa \cite{Cos}, de Figueiredo and Ding \cite{Fi-D}, Bartsch and Szulkin \cite{B-S} and Batkam and Colin \cite{B-C}.
\par In this paper we assume that $(\mathcal{S})$ is invariant under the action of the group $\mathcal{O}(N)$ of orthogonal transformation on $\mathbb{R}^N$, and we consider the existence of infinitely many radial and nonradial solutions. We recall that a function $u$ is radial if $u(x)=u(|x|)$ for every $x$. Inspired by Bartsch and willem \cite{BaWil}, we consider the restriction of the energy $J$ on suitable subspaces of invariant functions and we apply the generalized Fountain Theorem in \cite{B-C} to find infinitely many critical points of that restriction, which in turn are also critical points of $J$ by the Principle of Symmetric Criticality \cite{Pal}.
\par Before we state the main results of this paper, let us introduce our fundamental assumptions on the nonlinearity $F$:
\begin{itemize}
\item[$(F_1)$] $F\in \mathcal{C}^1(\mathbb{R}^N\times \mathbb{R}^2,\mathbb{R})$ and $F(x,0,0)=0$ $\forall x\in\mathbb{R}^N$.\\
\item[$(F_2)$] $|F_u(x,u,v)|+|F_v(x,u,v)|\leq c\big(|u|+|v|+|u|^{p-1}+|v|^{p-1}\big)$, with $2<p<2^*.$ \\
\item[$(F_3)$] $\exists\gamma>2$ such that $0<\gamma F(x,u,v)\leq uF_u(x,u,v)+vF_v(x,u,v),$ $\forall (u,v)\neq(0,0)$.\\
\item[$(F_4)$] $\inf\big\{F(x,u,v)\, \big| \, |(u,v)|\geq1,\, x\in\mathbb{R}^N\big\}>0$.\\
\item[$(F_5)$] $|F_u(x,u,v)|+|F_v(x,u,v)|=\circ(|(u,v)|)$, $|(u,v)|\rightarrow\infty$ uniformly on $\mathbb{R^N}$.\\
\item[$(F_6)$] $vF_v(x,u,v)\geq0$ or $uF_u(x,u,v)\geq0$ $\forall x\in\mathbb{R}^N$, $\forall(u,v)\in\mathbb{R}^2$.\\
\item[$(F_7)$] $F(x,u,v)=F(|x|,u,v)$, $\forall x\in\mathbb{R}^N$, $\forall(u,v)\in\mathbb{R}^2$.\\
\item[$(F_8)$] $F(x,-u,-v)=F(x,u,v)$, $\forall x\in\mathbb{R}^N$, $\forall(u,v)\in\mathbb{R}^2$.
\end{itemize}
Our main results are the following:
\begin{theo}\label{solutionradiale}
If $F$ satisfies the assumptions $(F_1)-(F_8)$, then $(\mathcal{S})$ has a sequence of radial solutions $(u_k,v_k)$ such that
\begin{equation*}
\frac{1}{2}\int_{\mathbb{R}^N}\big(|\nabla u_k|^2+|u_k|^2-|\nabla v_k|^2-|v_k|^2\big)dx-\int_{\mathbb{R}^N} F(x,u_k,v_k)dx\rightarrow+\infty, \quad k\rightarrow\infty.
\end{equation*}
\end{theo}

\begin{theo}\label{solutionnonradiale}
Let $N=4$ or $N\geq6$. If $F$ satisfied $(F_1)-(F_8)$, then $(\mathcal{S})$ has a sequence $(y_k,z_k)$ of nonradial solutions such that
\begin{equation*}
\frac{1}{2}\int_{\mathbb{R}^N}\big(|\nabla y_k|^2+|z_k|^2-|\nabla y_k|^2-|z_k|^2\big)dx-\int_{\mathbb{R}^N} F(x,y_k,z_k)dx\rightarrow+\infty, \quad k\rightarrow\infty.
\end{equation*}
\end{theo}
\par Similar results were obtained by  Huang and Li \cite{HuanLi} by applying the Limit Index Theory due to Li \cite{Li}. In this paper we use a different approach, which is much more direct and simpler, and which avoids approximation methods of Galerkin type. In particular, we do not need to verify a strong version of the Palais-Smale condition.
\par Our main results are proved in Section \ref{sectiontrois}, while the required abstract materials for the proofs are given is section \ref{sectiondeux}.

%%%%%%%%%%%%%%%%%%%%%%%%%%%%%%%%%%%%%%%%%%%%%%%%%%%%
%%%%%%%%%%%%%%%%%%%%%%%%%%%%%%%%%%%%%%%%%%%%%%%%%%%%
\section{Preliminaries}\label{sectiondeux}
%%%%%%%%%%%%%%%%%%%%%%%%%%%%%%%%%%%%%%%%%%%%%%%%%%%%
%%%%%%%%%%%%%%%%%%%%%%%%%%%%%%%%%%%%%%%%%%%%%%%%%%%%

In this section we recall the abstract results we will use in the proofs of the main theorems.
\subsection{Principle of symmetric criticality}
\begin{defin}
\textnormal{The action of a topological group $G$ on a normed vector space $X$ is a continuous map
\begin{equation*}
    G\times X\rightarrow X,\quad\quad (g,u)\mapsto g\cdot u
\end{equation*}
such that
\begin{equation*}
    1\cdot u=u, \,\, (gh)\cdot u=g(h\cdot u) \,\,\textnormal{and the map }u\mapsto g\cdot u \textnormal{ is linear, }\forall g,h\in G.
\end{equation*}
The action of $G$ is said to be isometric if $\|g\cdot u\|=\|u\|$, for every $u\in X$, $g\in G$.
A subset $A$ of $X$ is invariant if $g\cdot A=A$, for every $g\in G$. A function $\varphi:X\rightarrow \mathbb{R}$ is invariant if $\varphi(g\cdot u)=\varphi(u)$, for every $u\in X$, $g\in G$.}
\end{defin}
We denote
\begin{equation*}
    Fix(G):=\big\{u\in X\,\, ; \,\, g\cdot u=u\,\, \forall g\in G\big\}
\end{equation*}
the set of invariant points.\\
The following result is due to Palais \cite{Pal} (see also \cite{W}, Theorem $1.28$):
\begin{theo}[Principle of symmetric criticality]\label{criticality}
Assume that the action of the topological group $G$ on the Hilbert space $X$ is isometric. If $\varphi\in\mathcal{C}^1(X,\mathbb{R})$ is invariant and if $u$ is a critical point of $\varphi$ restricted to $Fix(G)$, then $u$ is a critical point of $\varphi$.
\end{theo}
\subsection{Generalized fountain theorem}
Now we assume that $X$ is a separable Hilbert space endowed with inner product $(\cdot)$  and the associated norm $\|\cdot\|$. Let $Y$ be a closed subspace of $X$. On $X=Y\oplus Y^\perp$ we consider the $\tau-$topology introduced by Kryszewski and Szulkin \cite{K-S}; that is the topology associated to the norm
\begin{equation*}
    \vvvert u\vvvert:=\max\big(\sum\limits_{j=0}^\infty|(Pu,a_j)|, \|Qu\|\big),\quad u\in X,
\end{equation*}
where $(a_j)_{j\geq0}$ is an orthonormal basis of $Y$, $P:X\rightarrow Y$ and $Q:X\rightarrow Z:=Y^\perp$ are the orthogonal projections. The $\tau-$topology has the following nice property \big(see \cite{K-S} or \cite{W}\big): If $(u_n)$ is a bounded sequence in $X$ then
\begin{equation*}
    u_n\stackrel{\tau}{\rightarrow}u\quad \Leftrightarrow \quad Pu_n\rightharpoonup Pu \quad \textnormal{and} \quad Qu_n\rightarrow Qu.
\end{equation*}
Consider an orthonormal basis $(e_j)_{j\geq0}$ of $Z$ and define
$$Y_k:=Y\oplus(\oplus_{j=0}^k{\mathbb R} e_j)\quad\quad \textnormal{ and }\quad\quad Z_k:=  \overline{\oplus_{j=k}^\infty\mathbb{R} e_j }. $$
\begin{defin}
\textnormal{ Let $\varphi\in\mathcal{C}^1(X,{\mathbb R})$.
\begin{enumerate}
              \item  $\varphi$ is said to satisfy the $(PS)_c$ condition \big(or the Palais-Smale condition at level $c$\big) if any sequence $(u_n)\subset X$ such that
\begin{equation*}
    \varphi(u_n)\rightarrow c \quad \textnormal{ and }\quad \varphi'(u_n)\rightarrow 0
\end{equation*}
has a convergent subsequence.
              \item We say that $\varphi$ is $\tau-$upper semicontinuous if the set $\big\{u\in X\,\, ;\,\, \varphi(u)\geq C\big\}$ is $\tau-$closed, for every $C\in\mathbb{R}$.
              \item We say that $\nabla\varphi$ is weakly sequentially continuous if the sequence $(\nabla\varphi (u_n))$ converges weakly to $\nabla\varphi (u)$ whenever $(u_n)$ converges weakly to $u$ in $X$.
            \end{enumerate}
}
\end{defin}
The following result is due to Batkam and Colin \cite{B-C}:
\begin{theo}[Generalized fountain theorem]\label{ft}
 Let $\varphi\in\mathcal{C}^1(X,{\mathbb R})$ be an even functional which is $\tau$-upper semicontinuous and such that $\nabla\varphi$ is weakly sequentially continuous. If there exist $\rho_k>r_k>0$ such that:
             \begin{itemize}
               \item [$(A_1)$]\,\,\,\  {$a_k \, := \, \displaystyle \sup_{\substack{u \in Y_k \\ \|u\| = \rho_k}} \varphi(u) \le 0 $ \,\,\ \textnormal{and} \,\,\ $\displaystyle \sup_{\substack{u \in Y_k \\ \|u\| \leq \rho_k}} \varphi(u)<\infty $}.
               \item [$(A_2)$] \,\,\,\ {$b_k \, := \, \displaystyle \inf_{\substack{u \in Z_k \\ \|u \| = r_k}} \varphi(u) \rightarrow \infty, \, k \rightarrow \infty$.}
               \item [$(A_3)$] \,\,\,\ $\varphi$ satisfies the $(PS)_c$ condition, $\forall c>0$.
             \end{itemize}
Then $\varphi$ has an unbounded sequence of critical values.
\end{theo}

%%%%%%%%%%%%%%%%%%%%%%%%%%%%%%%%%%%%%%%%%%%%%%%%%%%%
%%%%%%%%%%%%%%%%%%%%%%%%%%%%%%%%%%%%%%%%%%%%%%%%%%%%
\section{Proof of the main results}\label{sectiontrois}
%%%%%%%%%%%%%%%%%%%%%%%%%%%%%%%%%%%%%%%%%%%%%%%%%%%%
%%%%%%%%%%%%%%%%%%%%%%%%%%%%%%%%%%%%%%%%%%%%%%%%%%%%

Throughout this section $|\cdot|_p$ stands for the usual $L^p$ norm, and $\|\cdot\|$ stands for the usual $H^1$ norm. On $H^1\times H^1$ we consider the product norm $\|(u,v)\|^2=\|u\|^2+\|v\|^2.$ We denote $\rightarrow$ (resp. $\rightharpoonup$) the strong convergence (resp. the weak convergence).
\begin{defin}
Let $1\leq p,q<\infty$. On the space $L^p(\mathbb{R}^N)\cap L^q(\mathbb{R}^N)$ we define the norm
\begin{equation*}
    |u|_{p\wedge q}=|u|_p+|v|_q.
\end{equation*}
On the space $L^p(\mathbb{R}^N) + L^q(\mathbb{R}^N)$ we define the norm
\begin{equation*}
    |u|_{p\vee q}=\inf\big\{|v|_p+|w|_q;\,\,\, v\in L^p(\mathbb{R}^N),\, w\in L^q(\mathbb{R}^N),\,\,\, u=v+w\big\}.
\end{equation*}
\end{defin}
We refer to \cite{HuanLi} for the proof of the following lemma:
\begin{lem}\label{superposition}
Assume that $1\leq p,q, r, s<\infty$, $G\in \mathcal{C}(\mathbb{R}^N\times\mathbb{R}^2)$ and
\begin{equation*}
    |G(x,u,v)|\leq C\big(|u|^{\frac{p}{r}}+|v|^{\frac{p}{r}}+|u|^{\frac{q}{s}}+|v|^{\frac{q}{s}}\big).
\end{equation*}
Then, for every $u,v\in L^p(\mathbb{R}^N)\cap L^q(\mathbb{R}^N)$, $G(\cdot,u,v)\in L^r(\mathbb{R}^N) + L^s(\mathbb{R}^N)$ and the operator
\begin{equation*}
    T:\big(L^p(\mathbb{R}^N)\cap L^q(\mathbb{R}^N)\big)\times\big(L^p(\mathbb{R}^N)\cap L^q(\mathbb{R}^N)\big)\rightarrow L^r(\mathbb{R}^N) + L^s(\mathbb{R}^N)
\end{equation*}
defined by $T(u,v)=G(x,u,v)$ is continuous.
\end{lem}
\begin{lem}\label{regularite}
$J\in\mathcal{C}^1(H^1(\mathbb{R}^N)\times H^1(\mathbb{R}^N), \mathbb{R})$ with
\begin{equation}\label{jprime2}
\big<J'(u,v),(\phi,\varphi)\big>=(u,\phi)_1-(v,\varphi)_1-\int_{\mathbb{R}^N}\big(\phi F_u(x,u,v)+\varphi F_v(x,u,v)\big)dx,
\end{equation}
where $(\cdot,\cdot)_1$ denotes the usual inner product of $H^1(\mathbb{R}^N)$.
\end{lem}
\proof
\textbf{Existence of the Gateaux derivative.} Let $(u,v)\in H^1(\mathbb{R}^N)\times H^1(\mathbb{R}^N)$. For every $\phi,\varphi\in H^1(\mathbb{R}^N)$ and for $0<|t|<1$ we have:
\begin{multline*}
\frac{1}{t}\big[J(u+t\phi ,v+t\varphi)-J(u,v)\big]=\int_{\mathbb{R}^N}\big(\nabla u\nabla\phi+u\phi-\nabla v\nabla\varphi-v\varphi\big)dx\\
+\frac{t}{2}\int_{\mathbb{R}^N}\big(|\nabla\phi|^2+|\phi|^2-|\nabla\varphi|^2-|\varphi|^2\big)dx\\
-\int_{\mathbb{R}^N}\frac{1}{t}\big(F(x,u+t\phi,v+t\varphi)-F(x,u,v)\big)dx.
\end{multline*}
It follows from the mean value theorem that there exists $\lambda\in(0,1)$ such that, given $x\in \mathbb{R}^N$
\begin{multline*}
\frac{1}{|t|}|F(x,u+t\phi,v+t\varphi)-F(x,u,v)| \leq |F_u(x,u+t\lambda\phi,v+t\lambda\varphi)||\phi|\\+|F_v(x,u+t\lambda\phi,v+t\lambda\varphi)||\varphi|\\
\leq  \big(|F_u(x,u+t\lambda\phi,v+t\lambda\varphi)|+|F_v(x,u+t\lambda\phi,v+t\lambda\varphi)|\big)\big(|\phi|+|\varphi|\big)\\
\leq  c\big(|\lambda\phi|+|v+t\lambda\varphi|+|\lambda\phi|^{p-1}+|v+t\lambda\varphi|^{p-1}\big)\big(|\phi|+|\varphi|\big)\quad (\textnormal{in view of }(F_2))\\
\leq  c\big(|u|+|\phi|+|v|+|\varphi|+2^{p-1}(|u|^{p-1}+|\phi|^{p-1})+2^{p-1}(|v|^{p-1}+|\varphi|^{p-1})\big)\big(|\phi|+|\varphi|\big).
\end{multline*}
The H\"{o}lder inequality implies that the function $$c\big(|u|+|\phi|+|v|+|\varphi|+2^{p-1}(|u|^{p-1}+|\phi|^{p-1})+2^{p-1}(|v|^{p-1}+|\varphi|^{p-1})\big)(|\phi|+|\varphi|)$$
belongs to $L^1(\mathbb{R}^N).$ \\
It then follows from the dominated convergence theorem that
\begin{equation*}
\big<J'(u,v),(\phi,\varphi)\big>=(u,\phi)_1-(v,\varphi)_1-\int_{\mathbb{R}^N}\big(\phi F_u(x,u,v)+\varphi F_v(x,u,v)\big)dx.
\end{equation*}
\textbf{Continuity of the derivative.}
  Let $(u_n,v_n)\subset H^1(\mathbb{R}^N)\times H^1(\mathbb{R}^N)$ such that $(u_n,v_n)\rightarrow (u,v)$ in $H^1(\mathbb{R}^N)\times H^1(\mathbb{R}^N)$. By the Sobolev embedding theorem  $(u_n,v_n)\rightarrow (u,v)$ in $\big(L^2(\mathbb{R}^N\cap L^p(\mathbb{R}^N)\big)\times \big(L^2(\mathbb{R}^N\cap L^p(\mathbb{R}^N)\big)$. By Lemma \ref{superposition}  $F_u(x,u_n,v_n)\rightarrow F_u(x,u,v)$ and $F_v(x,u_n,v_n)\rightarrow F_v(x,u,v)$ in $L^2(\mathbb{R}^N)+L^{p'}(\mathbb{R}^N)$.\\
The H\"{o}lder inequality implies that
\begin{multline*}
|\big<J'(u_n,v_n)-J'(u,v),(\phi,\varphi)\big>|\leq\|u_n-u\|\|\phi\|+\|v_n-v\|\|\varphi\|\\
+|F_u(x,u_n,v_n)-F_u(x,u,v)|_{2\vee p'}|\phi|_{2\wedge p}\\
+|F_v(x,u_n,v_n)-F_v(x,u,v)|_{2\vee p'}|\varphi|_{2\wedge p}.
\end{multline*}
Hence we have
\begin{multline*}
\|J'(u_n,v_n)-J'(u,v)\|\leq \|u_n-u\|+\|v_n-v\|\\
+C\big[|F_u(x,u_n,v_n)-F_u(x,u,v)|_{2\vee p'}+|F_v(x,u_n,v_n)-F_v(x,u,v)|_{2\vee p'}\big].
\end{multline*}
We then deduce that $J'(u_n,v_n)-J'(u,v)\rightarrow0$ as $n\rightarrow\infty$.\qed

\subsection{Existence of radial solutions}

We recall that the action of the group $\mathcal{O}(N)$ on $H^1(\mathbb{R}^N)$ is defined by
\begin{equation*}
(g\cdot u)(x)=u(g^{-1}x).
\end{equation*}
 Let $H_{\mathcal{O}(N)}^1(\mathbb{R}^N):=\big\{u\in H^1(\mathbb{R}^N)\, \big|\, g\cdot u=u,\, \forall g\in \mathcal{O}(N)\big\}$ the set of radial functions. The following result is due to Strauss:
\begin{lem}[Strauss, 1977]\label{strauss}
Let $N\geq2$. For $2<q<2^\star$, the embedding $H_{\mathcal{O}(N)}^1(\mathbb{R}^N)\hookrightarrow L^q(\mathbb{R}^N)$ is compact.
\end{lem}
We define
\begin{equation*}
X:=H_{\mathcal{O}(N)}^1(\mathbb{R}^N)\times H_{\mathcal{O}(N)}^1(\mathbb{R}^N)\quad \textnormal{and} \quad \Phi:=J|_X.
\end{equation*}
 By $(F_7)$ $J$ is invariant under the action of $\mathcal{O}(N)$. It then follows from the Principle of Symmetric Criticality (Lemma \ref{criticality}) that the critical points of $\Phi$ are weak solutions of $(\mathcal{S})$.

\begin{lem}\label{palaissmale}
$\Phi$ satisfies the Palais-Smale condition on $X$. That is, every sequence $(u_n,v_n)\subset X$ such that $(\Phi(u_n,v_n))$ is bounded and $\Phi'(u_n,v_n)\rightarrow0$, has a convergent subsequence.
\end{lem}
\proof
Let $(u_n,v_n)\subset X$ such that
\begin{equation*}
d:=\sup_{\substack{n}}|\Phi(u_n)|<\infty \quad \textnormal{and}\quad \Phi'(u_n,v_n)\rightarrow0 \quad \textnormal{as}\quad n\rightarrow\infty.
\end{equation*}
We want to show that $(u_n,v_n)$ has a convergent subsequence.\\
By \eqref{jprime2} and $(F_6)$ we have
\begin{equation*}
\big<-\Phi'(u_n,v_n),(0,v_n)\big>=\|v_n\|^2+\int_{\mathbb{R}^N}v_nF_v(x,u_n,v_n)dx\geq \|v_n\|^2.
\end{equation*}
Hence for $n$ big enough we have $\|v_n\|^2\leq\|v_n\|$. This shows that $(v_n)$ is bounded.\\
On the other hand we deduce from \eqref{j2} and \eqref{jprime2} that
\begin{multline*}
\Phi(u_n,v_n)-\frac{1}{\gamma}\big<\Phi'(u_n,v_n),(u_n,v_n)\big>=\big(\frac{1}{2}-\frac{1}{\gamma}\big)\|u_n\|^2-\big(\frac{1}{2}-\frac{1}{\gamma}\big)\|v_n\|^2\\
+\int_{\mathbb{R}^N}\Big[\frac{1}{\gamma}\big(u_nF_u(x,u_n,v_n)+v_nF_v(x,u_n,v_n)\big)-F(x,u_n,v_n)\Big]dx\\
  \geq \big(\frac{1}{2}-\frac{1}{\gamma}\big)\|u_n\|^2-\big(\frac{1}{2}-\frac{1}{\gamma}\big)\|v_n\|^2 \quad \big(\textnormal{in view of }(F_3)\big).
\end{multline*}
Hence for $n$ big enough we have
\begin{equation*}
 \big(\frac{1}{2}-\frac{1}{\gamma}\big)\|u_n\|^2-\big(\frac{1}{2}-\frac{1}{\gamma}\big)\|v_n\|^2\leq \|(u_n,v_n)\|+d.
\end{equation*}
Since $(v_n)$ is bounded, we deduce that $(u_n)$ is also bounded.\\
\par Now up to a subsequence we have $(u_n,v_n)\rightarrow (u,v)$ in $X$. By Lemma \ref{strauss}, $u_n\rightarrow u$ and $v_n\rightarrow v$ in $L^P(\mathbb{R}^N)$.\\
We easily deduce from \eqref{j2} and \eqref{jprime2} that
\begin{multline*}
    \|u_n-u\|^2=\big<\Phi'(u_n,v_n)-\Phi'(u,v),(u_n-u,0)\big>\\
    +\int_{\mathbb{R}^N}\Big(F_u(x,u_n,v_n)-F_u(x,u,v)\Big)(u_n-u)dx,
\end{multline*}
\begin{multline*}
   \|v_n-v\|^2=-\big<\Phi'(u_n,v_n)-\Phi'(u,v),(0,v_n-v)\big>\\
   -\int_{\mathbb{R}^N}\Big(F_v(x,u_n,v_n)-F_v(x,u,v)\Big)(v_n-v)dx.
\end{multline*}
Clearly $\big<\Phi'(u_n,v_n)-\Phi'(u,v),(u_n-u,0)\big>\rightarrow0$ as $n\rightarrow\infty$.\\
$(F_2)$ and $(F_5)$ imply that for all $\varepsilon>0$ there exists $c_\varepsilon>0$ such that
\begin{equation*}
|F_u(x,u,v)|+|F_v(x,u,v)|\leq\varepsilon\big(|u|+|v|\big)+c_\varepsilon\big(|u|^{p-1}+|v|^{p-1}\big).
\end{equation*}
This implies that
\begin{multline*}
\big|\big(F_u(x,u_n,v_n)-F_u(x,u,v)\big)(u_n-u)\big|\leq \varepsilon \big(|u_n|+|v_n|+|u|+|v|\big)(|u_n|+|u|)\\
+c_\varepsilon\big(|u_n|^{p-1}+|u|^{p-1}+|v_n|^{p-1}+|v|^{p-1}\big)|u_n-u|.
\end{multline*}
Since the sequence $(u_n,v_n)$ is bounded in $X$, we obtain by using the H\"{o}lder inequality
\begin{equation*}
\int_{\mathbb{R}^N}\big|\big(F_u(x,u_n,v_n)-F_u(x,u,v)\big)(u_n-u)\big|dx\leq C\big(\varepsilon+c_\varepsilon|u_n-u|_p\big),
\end{equation*}
where $C$ is constant independent of $\varepsilon$ and $n$. It is then easy to see that $u_n\rightarrow u$ as $n\rightarrow \infty$. \\
By the same way we show that $v_n\rightarrow v$ as $n\rightarrow\infty$.\qed
\begin{lem}\label{faiblementsc}
$\nabla\Phi$ is weakly sequentially continuous.
\end{lem}
\proof
Let $(u_n,v_n)\subset X$ such that $(u_n,v_n)\rightharpoonup (u,v)$ in $X$.
 \eqref{jprime2} implies, for any $\phi,\varphi\in X$
 \begin{multline*}
    \big|\big<J'(u_n,v_n)-J'(u,v),(\phi,\varphi)\big>\big|\leq\Big|\int_{\mathbb{R}^N}\big(\nabla (u_n-u)\nabla\phi + (u_n-u)\phi\big)dx\Big|\\
    +\Big|\int_{\mathbb{R}^N}\big(\nabla (v_n-v)\nabla\varphi +(v_n-v)\varphi\big)dx\Big|
    + \Big|\int_{\mathbb{R}^N}\big(F_u(x,u_n,v_n)-F_u(x,u,v)\big)\phi dx\Big|\\
    +\Big|\int_{\mathbb{R}^N}\big(F_v(x,u_n,v_n)-F_v(x,u,v)\big)\varphi dx\Big|.
 \end{multline*}
It is clear that
\begin{equation*}
    \int_{\mathbb{R}^N}\big(\nabla (u_n-u)\nabla\phi + (u_n-u)\phi\big)dx\rightarrow0 \quad\textnormal{and}\quad \int_{\mathbb{R}^N}\big(\nabla (v_n-v)\nabla\varphi +(v_n-u)\varphi\big)dx\rightarrow0.
\end{equation*}
By Lemma \ref{strauss}, $u_n\rightarrow u$ and $v_n\rightarrow v$ in $L^p(\mathbb{R}^N)$. By $(F_2)$, $F_u$ and $F_v$ satisfy the assumption of Lemma \ref{superposition} with $r=q=p$ and $s=p':=\frac{p}{p-1}$. It then follows from Lemma \ref{superposition} that $F_u(x,u_n,v_n)-F_u(x,u,v)\rightarrow0$ and $F_v(x,u_n,v_n)-F_v(x,u,v)\rightarrow0$ in $L^p(\mathbb{R}^N) + L^{p'}(\mathbb{R}^N)$. Now, by using the H\"{o}lder's inequality we deduce that
\begin{equation*}
    \Big|\int_{\mathbb{R}^N}\big(F_u(x,u_n,v_n)-F_u(x,u,v)\big)\phi dx\Big|\leq\Big|F_u(x,u_n,v_n)-F_u(x,u,v)\Big|_{p\vee p'}|\phi|_{p\wedge p'}\rightarrow0
\end{equation*}
and
\begin{equation*}
    \Big|\int_{\mathbb{R}^N}\big(F_v(x,u_n,v_n)-F_v(x,u,v)\big)\varphi dx\Big|\leq\Big|F_v(x,u_n,v_n)-F_v(x,u,v)\Big|_{p\vee p'}|\varphi|_{p\wedge p'}\rightarrow0.
\end{equation*}
It then follows that $J'(u_n,v_n)\rightharpoonup J'(u,v)$.\qed

We define
\begin{equation*}
Y:=\{0\}\times H_{\mathcal{O}(N)}^1(\mathbb{R}^N), \quad \quad Z:= H_{\mathcal{O}(N)}^1(\mathbb{R}^N)\times\{0\},
\end{equation*}
and we consider the $\tau-$topology on $X=Y\oplus Z$.
\begin{lem}\label{tauscs}
$\Phi$ is $\tau-$upper semicontinuous.
\end{lem}
\proof
Let $(u_n,v_n)\subset X$ and $C\in\mathbb{R}$ such that $(u_n,v_n)\stackrel{\tau}{\rightarrow}(u,v)$ in $X$ and $J(u_n,v_n)\geq C$. By the definition of $\tau$, $u_n\rightarrow u$ in $H_{\mathcal{O}(N)}^1(\mathbb{R}^N)$.
\begin{equation*}
   J(u_n,v_n)\geq C \Longleftrightarrow\frac{1}{2}\|v_n\|^2+C\leq\frac{1}{2}\|u_n\|^2-\int_{\mathbb{R}^N} F(x,u_n,v_n)dx.
\end{equation*}
Since $F\geq0$ we deduce that
\begin{equation*}
    \frac{1}{2}\|v_n\|^2+C\leq\frac{1}{2}\|u_n\|^2.
\end{equation*}
Since $(u_n)$ is bounded, we easily deduce that $(v_n)$ is also bounded.\\
Now we may suppose, up to a subsequence that
\begin{eqnarray*}
    (u_n,v_n)&\rightharpoonup&(u,v)\, \textnormal{ in }\,X,\\
    u_n(x)&\rightarrow& u(x),\quad v_n(x)\rightarrow v(x)\,\,\, a.e. \, \textnormal{ in }\,\mathbb{R}^N, \\
    F(x,u_n(x),v_n(x))&\rightarrow& F(x,u(x),v(x))\,\,\, a.e. \, \textnormal{ in }\,\mathbb{R}^N.
\end{eqnarray*}
It follows from Fatou's lemma and the weak lower semicontinuity of the norm $\|\cdot\|$ that $-J(u,v)\leq -C$.\qed

\proof[\textbf{Theorem} \ref{solutionradiale}]
Let $(e_j)_{j\geq0}$ be an orthonormal basis of $H_{\mathcal{O}(N)}^1(\mathbb{R}^N)$. Set for $k\geq2$,
\begin{equation*}
    Y_k=Y\oplus\Big(\bigoplus\limits_{j=0}^k\mathbb{R}e_j\times\{0\}\Big) \,\,\, \textnormal{and} \,\,\ Z_k=\overline{\bigoplus\limits_{j=k}^\infty\mathbb{R}e_j}\times\{0\}.
\end{equation*}
 Let $(u,v)\in Y_k$. By $(F_3)$ and $(F_4)$, for all $\delta>0$ there exists $C_1=C_1(\delta)>0$ such that $F(x,u,v)\geq C_1|(u,v)|^\gamma-\delta|(u,v)|^2$.
 This implies that
 \begin{eqnarray*}
 % \nonumber to remove numbering (before each equation)
   \Phi(u,v) &\leq& \frac{1}{2}\|u\|^2-\frac{1}{2}\|v\|^2+\delta\big(|u|_2^2+|v|_2^2\big)-C_1|u|_\gamma^\gamma \\
             &\leq& \big(\frac{1}{2}+\delta\big)\|u\|^2+\big(\delta-\frac{1}{2}\big)\|v\|^2-C_1|u|_\gamma^\gamma.
 \end{eqnarray*}
Since all norms are equivalent on $\bigoplus\limits_{j=0}^k\mathbb{R}e_j$, we obtain
\begin{equation*}
\Phi(u,v)\leq \big(\frac{1}{2}+\delta\big)\|u\|^2+\big(\delta-\frac{1}{2}\big)\|v\|^2-C_2C_1\|u\|^\gamma,
\end{equation*}
where $C_2>0$ is a constant. By choosing $\delta<\frac{1}{4}$, we obtain
\begin{equation*}
\Phi(u,v)\leq \frac{3}{4}\|u\|^2-\frac{1}{4}\|v\|^2-C\|u\|^\gamma.
\end{equation*}
Hence $\Phi(u,v)\rightarrow-\infty$ as $\|(u,v)\|\rightarrow+\infty$, and assumption $(A_1)$ of Theorem \ref{ft} is then satisfied for $\rho_k$ sufficiently large.\\
\indent Now let $(u,0)\in Z_k$. Let $\varepsilon>0$, by $(F_2)$ and $F_3$ there exists $C_\varepsilon$ such that $F(x,u,0)\leq\varepsilon|u|^2+C_\varepsilon|u|^p$, which implies
\begin{equation*}
\Phi(u,0)\geq\frac{1}{2}\|u\|^2-\varepsilon|u|^2-C_\varepsilon |u|^p\geq \big(\frac{1}{2}-\varepsilon\big)\|u\|^2-C_\varepsilon |u|^p.
\end{equation*}
By choosing $\varepsilon<\frac{1}{4}$ we obtain
\begin{equation*}
\Phi(u,0)\geq \frac{1}{4}\|u\|^2-C |u|^p\geq \frac{1}{4}\|u\|^2-C \beta_k^p\|u\|^p,
\end{equation*}
where
\begin{equation*}
    \beta_k:=\sup_{\substack{w\in\overline{\oplus_{j=k}^\infty\mathbb{R} e^j}\\ \|w\|=1}}|w|_p.
\end{equation*}
Let
\begin{equation*}
    r_k:=\big(2pC\beta_k\big)^{\frac{1}{2-p}}.
\end{equation*}
Then for any $(u,0)\in Z_k$ such that $\|u\|=r_k$ we obtain
\begin{equation*}
    \Phi(u,0)\geq\big(\frac{1}{4}-\frac{1}{2p}\big)r_k^2\rightarrow\infty,\quad \textnormal{as}\quad k\rightarrow\infty,
\end{equation*}
since $\beta_k\rightarrow0$ as $k\rightarrow\infty$ \big(see \cite{W} Lemma $3.8$\big). The assumption $(A_2)$ of Theorem \ref{ft} is satisfied.\\
By $(F_8)$ $\Phi$ is even, and by Lemma \ref{palaissmale} the assumption $(A_3)$ of Theorem \ref{ft} is satisfied. We then conclude, in view of Lemmas \ref{regularite}, \ref{faiblementsc} and \ref{tauscs} by applying Theorem \ref{ft}.

\subsection{Existence of nonradial solutions}

Let $N=4$ or $N\geq6$. Let $2\leq m\leq N/2$ an integer different from $(N-1)/2$. We recall that the action of the group $G:=\mathcal{O}(m)\times\mathcal{O}(m)\times\mathcal{O}(N-2m)$ on $H^1(\mathbb{R}^N)$ is defined by $(g\cdot u)(x)=u(g^{-1}x)$.\\
 Let
\begin{equation*}
H_G^1(\mathbb{R}^N):=\big\{u\in H^1(\mathbb{R}^N)\,;\, g\cdot u=u \,\, \forall g\in G\big\}.
\end{equation*}
\begin{lem}[P. L. Lions, \cite{Lions}]
For $2<q<2^*$, the following embedding $H_G^1(\mathbb{R}^N)\hookrightarrow L^q(\mathbb{R}^N)$ is compact.
\end{lem}
Denote $\iota$ the involution defined on $\mathbb{R}^N=\mathbb{R}^m\times\mathbb{R}^m\times\mathbb{R}^{N-2m}$ by
$$\iota(x_1,x_2,x_3):=(x_2,x_1,x_3)$$
The action of $H:=\{id_{\mathbb{R}^N},\iota\}$ on $H_G^1(\mathbb{R}^N)$ is defined by
\begin{equation*}
(h\cdot u)(x)=\left\{
       \begin{array}{ll}
        u(x)\quad \textnormal{if }h=id_{\mathbb{R}^N}, & \hbox{} \\
        -u(h^{-1}x)\quad \textnormal{if }h=\iota. & \hbox{}
       \end{array}
     \right.
\end{equation*}
By this construction, $0$ is the only radial function in
\begin{equation*}
H_{G,H}^1(\mathbb{R}^N):=\big\{u\in H_G^1(\mathbb{R}^N)\,;\, h\cdot u=u,\,\forall h\in H\big\}.
\end{equation*}
\proof[\textbf{Theorem} \ref{solutionnonradiale}]
We define
\begin{equation*}
X:=H_{G,H}^1(\mathbb{R}^N)\times H_{G,H}^1(\mathbb{R}^N)\quad \textnormal{and } \Psi:=J|_{X}.
\end{equation*}
By $(F_7)$ and the Principe of Symmetric Criticality (Theorem \ref{criticality}), the critical points of $\Psi$ are also critical points of $J$.\\
Consider the $\tau-$topology on $X=Y\oplus Z$, where
\begin{equation*}
Y:=\{0\}\times H_{G,H}^1(\mathbb{R}^N) \textnormal{ and } Z:=H_{G,H}^1(\mathbb{R}^N)\times\{0\}.
\end{equation*}
The rest of the proof follows the lines of the proof of Theorem \ref{solutionradiale} above.\qed

\end{document}